\documentclass[11pt]{amsart}

\usepackage[margin=1in]{geometry}

\usepackage[utf8]{inputenc}
\usepackage[dvipsnames]{xcolor}
\usepackage[colorlinks=false,linkbordercolor={white},citebordercolor={white}]{hyperref}
\usepackage[backend=bibtex,backref=true,style=alphabetic,maxnames=4,maxalphanames=4]{biblatex}
\usepackage[capitalise,noabbrev]{cleveref}

\renewbibmacro{in:}{}
\usepackage{amscd}
\usepackage{amsmath}
\usepackage{amssymb}
\usepackage{comment}
\usepackage{enumitem}
\usepackage{mathrsfs}
\usepackage{mathtools}
\usepackage{relsize}
\usepackage{stmaryrd}
\usepackage{tikz-cd}
\usepackage{tikz}
\usepackage{url}
\usepackage{microtype}
\usepackage{array}

\usepackage{algorithm}
\usepackage{algpseudocode}
\algdef{SE}[SUBALG]{Indent}{EndIndent}{}{\algorithmicend\ }%
\algtext*{Indent}
\algtext*{EndIndent}
\numberwithin{algorithm}{section}

\numberwithin{equation}{section}
\newtheorem{lemma}[equation]{Lemma}
\newtheorem{theorem}[equation]{Theorem}
\newtheorem{proposition}[equation]{Proposition}
\newtheorem{corollary}[equation]{Corollary}

\newtheorem{claim*}{Claim}

\theoremstyle{definition}

\newtheorem{example}[equation]{Example}
\newtheorem{conventions}[equation]{Conventions}

\newtheorem{setup}[equation]{Setup}

\newtheorem{remark}[equation]{Remark}

\newcommand{\mfrak}[1]{\mathfrak{#1}}

\renewcommand{\k}{\Bbbk}

\newcommand{\m}{\mfrak{m}}

\newcommand{\Hom}{\operatorname{Hom}}

\newcommand{\D}{\msf{D}}
\newcommand{\kk}{\mathbf{k}}

\newcommand{\del}{\partial}
\newcommand{\Mod}{\operatorname{Mod}}

\DeclareMathOperator{\cone}{cone\!}
\DeclareMathOperator{\rank}{rank}

\def\nc{\newcommand}
\nc{\on}{\operatorname}
\nc{\bideg}{\on{bideg}}
\nc{\xra}{\xrightarrow}
\def\phi{\varphi}
\nc\cB{\mathcal{B}}
\def\th{{\on{th}}}

\def\D{\on{D}}
\def\Db{\D^{\on{b}}}

\nc{\into}{\hookrightarrow}
\nc{\onto}{\twoheadrightarrow}
\nc{\LL}{\mathbf{L}}
\nc{\RR}{\mathbf{R}}
\nc{\Perf}{\on{Perf}_{\on{gr}}}
\nc{\nat}{\natural}
\nc{\tors}{\on{tors}}
\nc{\Tors}{\on{Tors}}

\def\Mod{\on{Mod}}
\nc{\qgr}{\on{qgr}}
\nc{\Qgr}{\on{Qgr}}
\nc{\fQgr}{\on{Qgr}^{\on{f}}}
\nc{\colim}{\on{colim}}
\def\ZZ{\mathbb{Z}}
\nc{\Ext}{\on{Ext}}

\nc{\om}{\omega}
\nc{\w}{\widetilde}
\nc{\PP}{\mathbb{P}}
\nc{\mf}{\on{mf}}
\nc{\OO}{\mathcal{O}}
\nc{\Proj}{\on{Proj}}
\nc{\Qcoh}{\on{Qcoh}}
\nc{\coh}{\on{coh}}
\nc{\Tor}{\on{Tor}}
\nc{\Modf}{\Mod^{\on{f}}}

\nc{\ce}{\coloneqq}
\def\co{\colon}
\def\k{\kk}
\nc{\Com}{\on{Com}}
\nc{\A}{\mathcal{A}}
\nc{\B}{\mathcal{B}}
\nc{\CC}{\mathcal{C}}
\nc{\I}{\mathcal{I}}
\nc{\M}{\mathcal{M}}
\nc{\Sh}{\on{Sh}}
\nc{\QCoh}{\on{Qcoh}}
\nc{\Coh}{\on{coh}}
\nc{\fQCoh}{\QCoh^{\on{f}}}
\nc{\ov}{\overline}
\nc{\End}{\on{\underline{End}}}

\nc{\Qgrf}{\Qgr^{\on{f}}}
\nc{\uHom}{\underline{\Hom}}
\nc{\Inj}{\mathrm{Inj}}
\nc{\proj}{\mathrm{Proj}}
\nc{\spec}{\mathrm{Spec}}
\nc{\xla}{\xleftarrow}
\nc{\Dqgr}{\D_{\qgr}}
\nc{\DQgr}{\D_{\Qgr}}
\nc{\from}{\leftarrow}
\nc{\cd}{\on{cd}}
\nc{\N}{\mathcal{N}}
\nc{\F}{\mathcal{F}}
\nc{\FF}{\mathbb{F}}
\nc{\cC}{\mathcal{C}}
\nc{\cE}{\mathcal{E}}
\nc{\cF}{\mathcal{F}}
\nc{\cG}{\mathcal{G}}
\nc{\cH}{\mathcal{H}}
\nc{\cK}{\mathcal{K}}
\nc{\oa}{\overline{a}}
\nc{\ua}{\underline{a}}
\nc{\pd}{\on{pd}}
\nc{\wM}{\widetilde{M}}
\nc{\wN}{\widetilde{N}}
\nc{\wC}{\widetilde{C}}
\nc{\wD}{\widetilde{D}}
\nc{\wF}{\widetilde{F}}
\nc{\wG}{\widetilde{G}}
\nc{\cD}{\mathcal{D}}
\nc{\shHom}{\mathcal{H}om}
\nc{\ushHom}{\underline{\mathcal{H}om}}
\nc{\shExt}{\mathcal{E}xt}
\nc{\E}{\mathcal{E}}
\nc{\RHom}{\text{{\bf R}$\Hom$}}
\nc{\RuHom}{\text{{\bf R}$\underline{\Hom}$}}
\nc{\uExt}{\underline{\Ext}}
\nc{\cgeq}{\succcurlyeq}
\nc{\cle}{\prec}

\def\MR#1{}

\usepackage{listings}
\usepackage{xcolor}
\definecolor{maccolor}{rgb}{0.3,0.3,0.8}
\lstdefinelanguage{Macaulay2}
{
basicstyle={\ttfamily},
keywordstyle={\color{maccolor!80!black}},
commentstyle={\color{gray}},
stringstyle={\color{red!40!black}},
rulecolor=\color{maccolor},
basewidth={1.2ex}, 
sensitive=false,
morecomment=[l]{--},
morecomment=[s]{-*}{*-},
morestring=[b]",
escapechar={`},
escapebegin={\rmfamily},
morekeywords={about,abs,AbstractToricVarieties,accumulate,Acknowledgement,acos,acosh,acot,addCancelTask,addDependencyTask,addEndFunction,addHook,AdditionalPaths,addStartFunction,addStartTask,Adjacent,adjoint,AdjointIdeal,AffineVariety,AfterEval,AfterNoPrint,AfterPrint,agm,AInfinity,alarm,AlgebraicSplines,Algorithm,Alignment,all,AllCodimensions,allowableThreads,ambient,analyticSpread,Analyzer,AnalyzeSheafOnP1,ancestor,ancestors,ANCHOR,and,andP,AngleBarList,ann,annihilator,antipode,any,append,applicationDirectory,applicationDirectorySuffix,apply,applyKeys,applyPairs,applyTable,applyValues,apropos,argument,Array,arXiv,Ascending,ascii,asin,asinh,ass,assert,associatedGradedRing,associatedPrimes,AssociativeAlgebras,AssociativeExpression,atan,atan2,atEndOfFile,Authors,autoload,AuxiliaryFiles,backtrace,Bag,Bareiss,baseFilename,BaseFunction,baseName,baseRing,baseRings,BaseRow,BasicList,basis,BasisElementLimit,Bayer,BeforePrint,beginDocumentation,BeginningMacaulay2,Benchmark,benchmark,Bertini,BesselJ,BesselY,betti,BettiCharacters,BettiTally,between,BGG,BIBasis,Binary,BinaryOperation,Binomial,binomial,BinomialEdgeIdeals,Binomials,BKZ,BlockMatrix,BLOCKQUOTE,BODY,Body,BoijSoederberg,BOLD,Book3264Examples,Boolean,BooleanGB,borel,Boxes,BR,break,Browse,Bruns,cache,CacheExampleOutput,CacheFunction,CacheTable,cacheValue,CallLimit,cancelTask,capture,catch,Caveat,CC,CDATA,ceiling,Center,centerString,Certification,ChainComplex,chainComplex,ChainComplexExtras,ChainComplexMap,ChainComplexOperations,ChangeMatrix,char,CharacteristicClasses,characters,charAnalyzer,check,CheckDocumentation,chi,Chordal,class,Classic,clean,clearAll,clearEcho,clearOutput,close,closeIn,closeOut,ClosestFit,CODE,code,codim,CodimensionLimit,coefficient,CoefficientRing,coefficientRing,coefficients,Cofactor,CohenEngine,CohenTopLevel,CoherentSheaf,CohomCalg,cohomology,coimage,CoincidentRootLoci,coker,cokernel,collectGarbage,columnAdd,columnate,columnMult,columnPermute,columnRankProfile,columnSwap,combine,Command,commandInterpreter,commandLine,COMMENT,commonest,commonRing,comodule,CompactMatrix,compactMatrixForm,CompiledFunction,CompiledFunctionBody,CompiledFunctionClosure,Complement,complement,complete,CompleteIntersection,CompleteIntersectionResolutions,Complexes,ComplexField,components,compose,compositions,compress,concatenate,conductor,ConductorElement,cone,Configuration,ConformalBlocks,conjugate,connectionCount,Consequences,Constant,Constants,constParser,content,continue,contract,Contributors,ConvexInterface,conwayPolynomial,ConwayPolynomials,copy,copyDirectory,copyFile,copyright,Core,CorrespondenceScrolls,cos,cosh,cot,CotangentSchubert,cotangentSheaf,coth,cover,coverMap,cpuTime,createTask,Cremona,csc,csch,current,currentColumnNumber,currentDirectory,currentFileDirectory,currentFileName,currentLayout,currentLineNumber,currentPackage,currentString,currentTime,Cyclotomic,Database,Date,DD,dd,deadParser,debug,debugError,DebuggingMode,debuggingMode,debugLevel,DecomposableSparseSystems,Decompose,decompose,deepSplice,Default,default,defaultPrecision,Degree,degree,degreeLength,DegreeLift,DegreeLimit,DegreeMap,DegreeOrder,DegreeRank,Degrees,degrees,degreesMonoid,degreesRing,delete,demark,denominator,Dense,Density,Depth,depth,Descending,Descent,Describe,describe,Description,det,determinant,DeterminantalRepresentations,DGAlgebras,diagonalMatrix,diameter,Dictionary,dictionary,dictionaryPath,diff,DiffAlg,difference,dim,directSum,disassemble,discriminant,dismiss,Dispatch,distinguished,DIV,Divide,divideByVariable,DivideConquer,DividedPowers,Divisor,DL,Dmodules,do,doc,docExample,docTemplate,document,DocumentTag,Down,drop,DT,dual,eagonNorthcott,EagonResolution,echoOff,echoOn,EdgeIdeals,edit,EigenSolver,eigenvalues,eigenvectors,eint,EisenbudHunekeVasconcelos,elapsedTime,elapsedTiming,elements,Eliminate,eliminate,Elimination,EliminationMatrices,EllipticCurves,EllipticIntegrals,else,EM,Email,End,end,endl,endPackage,Engine,engineDebugLevel,EngineRing,EngineTests,entries,EnumerationCurves,environment,Equation,EquivariantGB,erase,erf,erfc,error,errorDepth,euler,EulerConstant,eulers,even,EXAMPLE,ExampleFiles,ExampleItem,examples,ExampleSystems,Exclude,exec,exit,exp,expectedReesIdeal,expm1,exponents,export,exportFrom,exportMutable,Expression,expression,Ext,extend,ExteriorIdeals,ExteriorModules,exteriorPower,Factor,factor,false,Fano,FastMinors,FastNonminimal,FGLM,File,fileDictionaries,fileExecutable,fileExists,fileExitHooks,fileLength,fileMode,FileName,FilePosition,fileReadable,fileTime,fileWritable,fillMatrix,findFiles,findHeft,FindOne,findProgram,findSynonyms,FiniteFittingIdeals,First,first,firstkey,FirstPackage,fittingIdeal,flagLookup,FlatMonoid,flatten,flattenRing,Flexible,flip,floor,flush,fold,FollowLinks,for,forceGB,fork,FormalGroupLaws,Format,format,formation,FourierMotzkin,FourTiTwo,fpLLL,frac,fraction,FractionField,frames,freeResolution,FrobeniusThresholds,from,fromDividedPowers,fromDual,Function,FunctionApplication,FunctionBody,functionBody,FunctionClosure,FunctionFieldDesingularization,fusePairs,futureParser,GaloisField,Gamma,gb,GBDegrees,gbRemove,gbSnapshot,gbTrace,gcd,gcdCoefficients,gcdLLL,GCstats,genera,GeneralOrderedMonoid,GenerateAssertions,generateAssertions,generator,generators,Generic,GenericInitialIdeal,genericMatrix,genericSkewMatrix,genericSymmetricMatrix,gens,genus,get,getc,getChangeMatrix,getenv,getGlobalSymbol,getNetFile,getNonUnit,getPrimeWithRootOfUnity,getSymbol,getWWW,GF,gfanInterface,Givens,GKMVarieties,GLex,Global,global,globalAssign,globalAssignFunction,GlobalAssignHook,globalAssignment,globalAssignmentHooks,GlobalDictionary,GlobalHookStore,globalReleaseFunction,GlobalReleaseHook,Gorenstein,GradedLieAlgebras,GradedModule,gradedModule,GradedModuleMap,gradedModuleMap,gramm,GraphicalModels,GraphicalModelsMLE,Graphics,graphIdeal,graphRing,Graphs,Grassmannian,GRevLex,GroebnerBasis,groebnerBasis,GroebnerBasisOptions,GroebnerStrata,GroebnerWalk,groupID,GroupLex,GroupRevLex,GTZ,Hadamard,handleInterrupts,HardDegreeLimit,hash,HashTable,hashTable,HEAD,HEADER1,HEADER2,HEADER3,HEADER4,HEADER5,HEADER6,HeaderType,Heading,Headline,Heft,heft,Height,height,help,Hermite,hermite,Hermitian,HH,hh,HigherCIOperators,HighestWeights,Hilbert,hilbertFunction,hilbertPolynomial,hilbertSeries,HodgeIntegrals,hold,Holder,Hom,homeDirectory,HomePage,Homogeneous,Homogeneous2,homogenize,homology,homomorphism,HomotopyLieAlgebra,hooks,horizontalJoin,HorizontalSpace,HR,HREF,HTML,html,httpHeaders,Hybrid,HyperplaneArrangements,Hypertext,hypertext,HypertextContainer,HypertextParagraph,icFracP,icFractions,icMap,icPIdeal,id,Ideal,ideal,idealizer,identity,if,IgnoreExampleErrors,ii,image,imaginaryPart,IMG,ImmutableType,importFrom,in,incomparable,Increment,independentSets,indeterminate,IndeterminateNumber,Index,index,indexComponents,IndexedVariable,IndexedVariableTable,indices,inducedMap,inducesWellDefinedMap,InexactField,InexactFieldFamily,InexactNumber,InfiniteNumber,infinity,info,InfoDirSection,infoHelp,Inhomogeneous,input,Inputs,insert,installAssignmentMethod,installedPackages,installHilbertFunction,installMethod,installMinprimes,installPackage,InstallPrefix,instance,instances,IntegralClosure,integralClosure,integrate,IntermediateMarkUpType,interpreterDepth,intersect,intersectInP,Intersection,intersection,interval,InvariantRing,inverse,InverseMethod,inversePermutation,Inverses,inverseSystem,InverseSystems,Invertible,InvolutiveBases,irreducibleCharacteristicSeries,irreducibleDecomposition,isAffineRing,isANumber,isBorel,isCanceled,isCommutative,isConstant,isDirectory,isDirectSum,isEmpty,isField,isFinite,isFinitePrimeField,isFreeModule,isGlobalSymbol,isHomogeneous,isIdeal,isInfinite,isInjective,isInputFile,isIsomorphism,isLinearType,isListener,isLLL,isMember,isModule,isMonomialIdeal,isNormal,isOpen,isOutputFile,isPolynomialRing,isPrimary,isPrime,isPrimitive,isPseudoprime,isQuotientModule,isQuotientOf,isQuotientRing,isReady,isReal,isReduction,isRegularFile,isRing,isSkewCommutative,isSorted,isSquareFree,isStandardGradedPolynomialRing,isSubmodule,isSubquotient,isSubset,isSupportedInZeroLocus,isSurjective,isTable,isUnit,isWellDefined,isWeylAlgebra,ITALIC,Iterate,Jacobian,jacobian,jacobianDual,Jets,Join,join,Jupyter,K3Carpets,K3Surfaces,Keep,KeepFiles,KeepZeroes,ker,kernel,kernelLLL,kernelOfLocalization,Key,keys,Keyword,Keywords,kill,koszul,Kronecker,KustinMiller,LABEL,last,lastMatch,LATER,LatticePolytopes,Layout,lcm,leadCoefficient,leadComponent,leadMonomial,leadTerm,Left,left,length,LengthLimit,letterParser,Lex,LexIdeals,LI,Licenses,LieTypes,lift,liftable,Limit,limitFiles,limitProcesses,Linear,LinearAlgebra,LinearTruncations,lineNumber,lines,LINK,linkFile,List,list,listForm,listLocalSymbols,listSymbols,listUserSymbols,LITERAL,LLL,LLLBases,lngamma,load,loadDepth,LoadDocumentation,loadedFiles,loadedPackages,loadPackage,Local,local,localDictionaries,LocalDictionary,localize,LocalRings,locate,log,log1p,LongPolynomial,lookup,lookupCount,LowerBound,LUdecomposition,M0nbar,M2CODE,Macaulay2Doc,makeDirectory,MakeDocumentation,makeDocumentTag,MakeHTML,MakeInfo,MakeLinks,makePackageIndex,MakePDF,makeS2,Manipulator,map,MapExpression,MapleInterface,markedGB,Markov,MarkUpType,match,mathML,Matrix,matrix,MatrixExpression,Matroids,max,maxAllowableThreads,maxExponent,MaximalRank,maxPosition,MaxReductionCount,MCMApproximations,member,memoize,memoizeClear,memoizeValues,MENU,merge,mergePairs,META,method,MethodFunction,MethodFunctionBinary,MethodFunctionSingle,MethodFunctionWithOptions,methodOptions,methods,midpoint,min,minExponent,mingens,mingle,minimalBetti,MinimalGenerators,MinimalMatrix,minimalPresentation,minimalPresentationMap,minimalPresentationMapInv,MinimalPrimes,minimalPrimes,minimalReduction,Minimize,minimizeFilename,MinimumVersion,minors,minPosition,minPres,minprimes,Minus,minus,Miura,MixedMultiplicity,mkdir,mod,Module,module,ModuleDeformations,modulo,MonodromySolver,Monoid,monoid,MonoidElement,Monomial,MonomialAlgebras,monomialCurveIdeal,MonomialIdeal,monomialIdeal,MonomialIntegerPrograms,MonomialOrbits,MonomialOrder,Monomials,monomials,MonomialSize,monomialSubideal,moveFile,multidegree,multidoc,multigraded,MultigradedBettiTally,MultiGradedRationalMap,multiplicity,MultiplicitySequence,MultiplierIdeals,MultiplierIdealsDim2,MultiprojectiveVarieties,mutable,MutableHashTable,mutableIdentity,MutableList,MutableMatrix,mutableMatrix,NAGtypes,Name,nanosleep,Nauty,NautyGraphs,NCAlgebra,NCLex,needs,needsPackage,Net,net,NetFile,netList,new,newClass,newCoordinateSystem,NewFromMethod,newline,NewMethod,newNetFile,NewOfFromMethod,NewOfMethod,newPackage,newRing,nextkey,nextPrime,nil,NNParser,NoetherianOperators,NoetherNormalization,NonAssociativeProduct,NonminimalComplexes,nonspaceAnalyzer,NoPrint,norm,normalCone,Normaliz,NormalToricVarieties,not,Nothing,notify,notImplemented,NTL,null,nullaryMethods,nullhomotopy,nullParser,nullSpace,Number,number,NumberedVerticalList,numcols,numColumns,numerator,numeric,NumericalAlgebraicGeometry,NumericalCertification,NumericalImplicitization,NumericalLinearAlgebra,NumericalSchubertCalculus,numericInterval,NumericSolutions,numgens,numRows,numrows,odd,oeis,of,ofClass,OL,OldPolyhedra,OldToricVectorBundles,on,OneExpression,OnlineLookup,OO,oo,ooo,oooo,openDatabase,openDatabaseOut,openFiles,openIn,openInOut,openListener,OpenMath,openOut,openOutAppend,operatorAttributes,Option,OptionalComponentsPresent,optionalSignParser,Options,options,OptionTable,optP,or,Order,order,OrderedMonoid,orP,OutputDictionary,Outputs,override,pack,Package,package,PackageCitations,PackageDictionary,PackageExports,PackageImports,PackageTemplate,packageTemplate,pad,pager,PairLimit,pairs,PairsRemaining,PARA,Parametrization,parent,Parenthesize,Parser,Parsing,part,Partition,partition,partitions,parts,path,pdim,peek,PencilsOfQuadrics,Permanents,permanents,permutations,pfaffians,PHCpack,PhylogeneticTrees,pi,PieriMaps,pivots,PlaneCurveSingularities,plus,poincare,poincareN,Points,polarize,poly,Polyhedra,Polymake,PolynomialRing,Posets,Position,position,positions,PositivityToricBundles,POSIX,Postfix,Power,power,powermod,PRE,Precision,precision,Prefix,prefixDirectory,prefixPath,preimage,prepend,presentation,pretty,primaryComponent,PrimaryDecomposition,primaryDecomposition,PrimaryTag,PrimitiveElement,Print,print,printerr,printingAccuracy,printingLeadLimit,printingPrecision,printingSeparator,printingTimeLimit,printingTrailLimit,printString,printWidth,processID,Product,product,ProductOrder,profile,profileSummary,Program,programPaths,ProgramRun,Proj,Projective,ProjectiveHilbertPolynomial,projectiveHilbertPolynomial,ProjectiveVariety,promote,protect,Prune,prune,PruneComplex,pruningMap,Pseudocode,pseudocode,pseudoRemainder,Pullback,PushForward,pushForward,Python,QQ,QQParser,QRDecomposition,QthPower,Quasidegrees,QuaternaryQuartics,QuillenSuslin,quit,Quotient,quotient,quotientRemainder,QuotientRing,Radical,radical,RadicalCodim1,radicalContainment,RaiseError,random,RandomCanonicalCurves,RandomComplexes,RandomCurves,RandomCurvesOverVerySmallFiniteFields,RandomGenus14Curves,RandomIdeals,randomKRationalPoint,RandomMonomialIdeals,randomMutableMatrix,RandomObjects,RandomPlaneCurves,RandomPoints,RandomSpaceCurves,Range,rank,RationalMaps,RationalPoints,RationalPoints2,ReactionNetworks,read,readDirectory,readlink,readPackage,RealField,RealFP,realPart,realpath,RealQP,RealQP1,RealRoots,RealRR,RealXD,recursionDepth,recursionLimit,Reduce,reducedRowEchelonForm,reduceHilbert,reductionNumber,ReesAlgebra,reesAlgebra,reesAlgebraIdeal,reesIdeal,References,ReflexivePolytopesDB,regex,regexQuote,registerFinalizer,regSeqInIdeal,Regularity,regularity,relations,RelativeCanonicalResolution,relativizeFilename,Reload,remainder,RemakeAllDocumentation,remove,removeDirectory,removeFile,removeLowestDimension,reorganize,replace,RerunExamples,res,reshape,ResidualIntersections,ResLengthThree,Resolution,resolution,ResolutionsOfStanleyReisnerRings,restart,Result,resultant,Resultants,return,returnCode,Reverse,reverse,RevLex,Right,right,Ring,ring,RingElement,RingFamily,ringFromFractions,RingMap,rootPath,roots,rootURI,rotate,round,rowAdd,RowExpression,rowMult,rowPermute,rowRankProfile,rowSwap,RR,RRi,rsort,run,RunDirectory,RunExamples,RunExternalM2,runHooks,runLengthEncode,runProgram,same,saturate,Saturation,scan,scanKeys,scanLines,scanPairs,scanValues,schedule,schreyerOrder,Schubert,Schubert2,SchurComplexes,SchurFunctors,SchurRings,SCRIPT,scriptCommandLine,ScriptedFunctor,SCSCP,searchPath,sec,sech,SectionRing,SeeAlso,seeParsing,SegreClasses,select,selectInSubring,selectVariables,SelfInitializingType,SemidefiniteProgramming,Seminormalization,separate,SeparateExec,separateRegexp,Sequence,sequence,Serialization,serialNumber,Set,set,setEcho,setGroupID,setIOExclusive,setIOSynchronized,setIOUnSynchronized,setRandomSeed,setup,setupEmacs,sheaf,SheafExpression,sheafExt,sheafHom,SheafOfRings,shield,ShimoyamaYokoyama,short,show,showClassStructure,showHtml,showStructure,showTex,showUserStructure,SimpleDoc,simpleDocFrob,SimplicialComplexes,SimplicialDecomposability,SimplicialPosets,SimplifyFractions,sin,singularLocus,sinh,size,size2,SizeLimit,SkewCommutative,SlackIdeals,sleep,SLnEquivariantMatrices,SLPexpressions,SMALL,smithNormalForm,solve,someTerms,Sort,sort,sortColumns,SortStrategy,source,SourceCode,SourceRing,SPACE,SpaceCurves,SPAN,span,SparseMonomialVectorExpression,SparseResultants,SparseVectorExpression,Spec,SpechtModule,SpecialFanoFourfolds,specialFiber,specialFiberIdeal,SpectralSequences,splice,splitWWW,sqrt,SRdeformations,stack,stacksProject,Standard,standardForm,standardPairs,StartWithOneMinor,stashValue,StatePolytope,StatGraphs,status,stderr,stdio,step,StopBeforeComputation,stopIfError,StopWithMinimalGenerators,Strategy,String,STRONG,StronglyStableIdeals,STYLE,Style,style,SUB,sub,SubalgebraBases,sublists,submatrix,submatrixByDegrees,Subnodes,subquotient,SubringLimit,Subscript,subscript,SUBSECTION,subsets,substitute,substring,subtable,Sugarless,Sum,sum,SumOfTwists,SumsOfSquares,SUP,super,SuperLinearAlgebra,Superscript,superscript,support,SVD,SVDComplexes,switch,SwitchingFields,sylvesterMatrix,Symbol,symbol,SymbolBody,symbolBody,SymbolicPowers,symlinkDirectory,symlinkFile,symmetricAlgebra,symmetricAlgebraIdeal,symmetricKernel,SymmetricPolynomials,symmetricPower,synonym,SYNOPSIS,syz,Syzygies,SyzygyLimit,SyzygyMatrix,SyzygyRows,syzygyScheme,TABLE,Table,table,take,Tally,tally,tan,TangentCone,tangentCone,tangentSheaf,tanh,target,Task,taskResult,TateOnProducts,TD,temporaryFileName,tensor,tensorAssociativity,TensorComplexes,terminalParser,terms,TEST,Test,testExample,testHunekeQuestion,TestIdeals,TestInput,tests,TEX,tex,TeXmacs,texMath,Text,TH,then,Thing,ThinSincereQuivers,ThreadedGB,threadVariable,Threshold,throw,Time,time,times,timing,TITLE,TO,to,TO2,toAbsolutePath,toCC,toDividedPowers,toDual,toExternalString,toField,TOH,toList,toLower,top,top,topCoefficients,Topcom,topComponents,topLevelMode,Tor,TorAlgebra,Toric,ToricInvariants,ToricTopology,ToricVectorBundles,toRR,toRRi,toSequence,toString,TotalPairs,toUpper,TR,trace,transpose,TriangularSets,Tries,Trim,trim,Triplets,Tropical,true,Truncate,truncate,truncateOutput,Truncations,try,TSpreadIdeals,TT,tutorial,Type,TypicalValue,typicalValues,UL,ultimate,unbag,uncurry,Undo,undocumented,uniform,uninstallAllPackages,uninstallPackage,Unique,unique,Units,Unmixed,unsequence,unstack,Up,UpdateOnly,UpperTriangular,URL,urlEncode,Usage,use,UseCachedExampleOutput,UseHilbertFunction,UserMode,userSymbols,UseSyzygies,utf8,utf8check,validate,value,values,Variable,VariableBaseName,Variables,Variety,variety,vars,Vasconcelos,Vector,vector,VectorExpression,VectorFields,VectorGraphics,Verbose,Verbosity,Verify,VersalDeformations,versalEmbedding,Version,version,VerticalList,VerticalSpace,viewHelp,VirtualResolutions,VirtualTally,VisibleList,Visualize,wait,WebApp,wedgeProduct,weightRange,Weights,WeylAlgebra,WeylGroups,when,whichGm,while,width,wikipedia,Wrap,wrap,WrapperType,XML,xor,youngest,zero,ZeroExpression,zeta,ZZ,ZZParser}
}
\lstalias{Macaulay2output}{Macaulay2}

\author{Michael K. Brown}
\author{Souvik Dey}
\author{Guanyu Li}
\author{Mahrud Sayrafi}

\newcommand{\Addresses}{
  \footnotesize
  \vskip\baselineskip
  \noindent \textsc{Department of Mathematics and Statistics, Auburn University}   \par\nopagebreak
  \noindent \textit{E-mail address:} \texttt{mkb0096@auburn.edu}
  \vskip\baselineskip
  \noindent \textsc{Department of Mathematical Sciences, University of Arkansas} \par\nopagebreak
  \noindent \textit{E-mail address:} \texttt{souvikd@uark.edu}
  \vskip\baselineskip
  \noindent \textsc{Department of Mathematics, Cornell University} \par\nopagebreak
  \noindent \textit{E-mail address:} \texttt{gl479@cornell.edu}
  \vskip\baselineskip
  \noindent \textsc{Department of Mathematics and Statistics, McMaster University} \par\nopagebreak
  \noindent \textit{E-mail address:} \texttt{mahrud@mcmaster.ca}
  \vskip\baselineskip
}

\makeatletter
\@namedef{subjclassname@2020}{%
  \textup{2020} Mathematics Subject Classification}
\makeatother

\begin{document}


\title{Computing global Ext for complexes}

\begin{abstract}
  We give a computational algorithm for computing Ext groups between bounded complexes of coherent sheaves on a projective variety, and we describe an implementation of this algorithm in \verb|Macaulay2|. In particular, our results yield methods for computing derived global sections of bounded complexes of coherent sheaves and mutations of exceptional collections. 
\end{abstract}

\numberwithin{equation}{section}

\maketitle

\setcounter{section}{0}


\vspace{-1.5em} 

\section{Introduction}\label{sec:intro}

Consider bounded complexes of coherent sheaves $\cC$ and $\cD$ on a projective variety $X \subseteq \PP^n$\!.\linebreak Our goal is to provide an algorithm for effectively computing the extension groups $\Ext_X^m(\cC, \cD)$ for~$m\in\ZZ$. If $\F$ and $\cG$ are coherent sheaves on $X$, and $F$ and $G$ are graded modules over the homogeneous coordinate ring $R$ of $X$ lifting $\F$ and $\cG$, an algorithm due to Smith~\cite{Smith00} computes $\Ext_X^m(\F,\cG)$ via an isomorphism with $\Ext_R^m(F_{\geq d},G)$ for high enough degree $d\in\ZZ$. We extend this algorithm to bounded complexes of coherent sheaves. 

Before stating our main theorem, we fix some notation that will be used throughout the paper. Let $\k$ be a field and $S = \k[x_0, \dots, x_n]$ be the $\ZZ$-graded polynomial ring with $\deg(x_i) = 1$. Given a graded $S$-module $M$\!, let $\pd_S(M)$ denote its projective dimension. 
Let $I \subseteq S$ be a homogeneous ideal, $R = S/I$, and $X = \Proj(R) \subseteq \PP^n$\!. Given a bounded complex $C$ of finitely generated graded $R$-modules, we write $C^j$ for the $j^\th$ cohomological term of~$C$, and we use $C_i \ce \bigoplus_{ j \in \ZZ} C^j_i$ to denote the complex of $\k$-vector spaces given by the $i^\th$ graded components of each cohomological term of~$C$. Further, we use $C_{\ge r}$ to denote the subcomplex  $\bigoplus_{i \ge r} C_i$ of $C$. We also write
\[
\inf(C) \ce \inf\{j \text{ : } C^j \ne 0\} \quad \text{and} \quad
\sup(C) \ce \sup\{j \text{ : } C^j \ne 0\}.
\]
We define the \emph{dimension of $C$} to be $\dim(\bigoplus_{j \in \ZZ} C^j)$. The Betti numbers of $C$ over $R$ are given by $\beta^R_{i,j}(C) \ce \dim_\k H^{-i}(C \otimes^\LL_R \k)_j$, and following \cite{Smith00} we define:
\[
\ua^R_i(C) \ce \inf\{ j \colon \beta^R_{i, j}(C) \ne 0\}, \quad
\oa^R_i(C) \ce \sup\{ j \colon \beta^R_{i, j}(C) \ne 0\}.
\]
The values $\oa^S_i(C)$ and $\ua^S_i(C)$ are defined in the same way, but involving Betti numbers over $S$.

The following is our main result, which is an analogue of \cite[Theorem 1]{Smith00} for complexes. 

\begin{theorem}\label{thm:smithcpxintro}
  Let $C$ and $D$ be bounded complexes of finitely generated $R$-modules such that $\sup(C) = \sup(D) = 0$. Let $m \in \ZZ$, and set $\ell \ce  \min\{\dim(D), m - \inf(D)\}$. If the inequality
  \[ r \ge \max\{\oa^S_{i}(D^{j}) \text{ $\on{:}$ } n-\ell \le i \le \pd_S(D^j), \text{ } \inf(D) \le j \le  0\} - n \]
  holds, then there is a canonical isomorphism of graded $R$-modules
  \[ \uExt^m_R(C_{\ge r}, D)_{\ge 0} \xra{\cong} \bigoplus_{v \ge 0} \Ext^m_X(\wC, \wD(v)). \]
  If we have the stronger inequality
 $r \ge \max\{\oa^S_{i}(D^{j}) \text{ $\on{:}$ } 0 \le i \le \pd_S(D^j), \text{ } \inf(D) \le j \le  0\} - n$, then there is a canonical quasi-isomorphism of complexes of graded $R$-modules
  \[ \RuHom_R(C_{\ge r}, D)_{\ge 0} \xra{\simeq} \bigoplus_{v \ge 0} \RHom_X(\wC, \wD(v)). \]
\end{theorem}

We refer to Conventions~\ref{notation} for the meaning of $\uExt_R$, $\RuHom_R$, and $\RHom_X$ in Theorem~\ref{thm:smithcpxintro}. Theorem~\ref{thm:smithcpxintro} does not quite recover Smith's theorem, because of the presence of $-n$ rather than~$-m$ in the inequalities. However, when $H^j(C) = 0$ for $j \ne 0$, Theorem~\ref{thm:smithcpxintro} can be strengthened to a result that does generalize~\cite[Theorem 1]{Smith00}: see Corollary~\ref{cor:positive} for the detailed statement. The condition $\sup(C) = \sup(D)  = 0$ is included only for simplicity; see Theorem~\ref{thm:smithcpx} for a more general statement without this assumption. 

We prove Theorem~\ref{thm:smithcpxintro} in Section~\ref{sec:proof}. The argument requires careful bookkeeping to sharpen the various bounds as much as possible. For instance, one of our key technical results, Proposition~\ref{prop:SheafExtToHom}, is proven by showing that certain terms of a hypercohomology spectral sequence stabilize at the second page. Theorem~\ref{thm:smithcpxintro} leads to an effective algorithm for computing $\Ext$ between bounded complexes of coherent sheaves on $X$\!, and in particular for computing their derived global sections. We have implemented this algorithm in the symbolic algebra system \verb|Macaulay2|~\cite{M2}, which we demonstrate through examples in Section~\ref{sec:examples}.

Finally, in Section~\ref{sec:exc-cols} we highlight applications of our results in the study of derived categories in algebraic geometry. For instance, one may use Theorem~\ref{thm:smithcpxintro} to check computationally whether a family of complexes of sheaves forms an exceptional collection in the derived category and to compute mutations of exceptional collections and spherical twists of complexes. In future work, we will apply our results to give an implementation in \verb|Macaulay2| of the fully faithful embeddings in Orlov's Landau--Ginzburg/Calabi--Yau correspondence~\cite[Theorem 2.5]{Orlov09}; indeed, this was our original motivation for the present work.

\subsection*{Acknowledgments}
We thank Devlin Mallory, Greg Smith, Keller VandeBogert, and Sasha Zotine for helpful conversations, and we thank the University of Utah for hosting the 2024 conference ``Computational Algebraic Geometry and String Theory'' (funded by NSF grant DMS-2001206), where work on this project began. The first author was partially supported by NSF grants DMS-2302373, DMS-2302375, and DMS-2412042. The second author was partially supported by the Charles University Research Center program No.~UNCE/24/SCI/022 and grant GA \v{C}R 23-05148S from the Czech Science Foundation. We also thank the referee for many helpful suggestions.

\begin{conventions}\label{notation}
  As indicated above, we index our complexes cohomologically. Let $C$ and $D$ be complexes of graded modules over a $\ZZ$-graded ring $A$ with differentials~$\del_C$ and $\del_D$. The \emph{$i^{\th}$ cohomological shift of $C$} is the complex $C[i]$ with $C[i]^j = C^{i+j}$ and differential $(-1)^i\del_C$. The \emph{$j^{\th}$ twist of~$C$} is the complex $C(j)$ with $C(j)_i \ce C_{i + j}$ and the same differential as $C$. We denote by $\Hom_A(C, D)$ the set of morphisms of complexes from $C$ to $D$ of internal and homological degree 0, and $\uHom_A(C, D)$ denotes the complex of graded $A$-modules whose $j^{\th}$ term is the graded $A$-module $\bigoplus_{i \in \ZZ} \Hom_A(C, D(i)[j])$ and whose differential sends a map~$g$ of cohomological degree $j$ to the map $\del_D g - (-1)^{j}g\del_C$. We write $\RuHom_A(C, D)$ for the derived $\Hom$ complex from $C$ to $D$, and we set
  \[ \uExt^m_A(C, D) \ce H^m \RuHom_A(C, D) \quad \text{and} \quad \Ext^m_A(C, D) \ce \uExt^m_A(C, D)_0. \]
Suppose now that $C$ and~$D$ are complexes of graded $R$-modules. We write $\RHom_X(\wC,\wD)$ for the complex of $\k$-vector spaces given by derived $\Hom$ and define the graded $R$-modules
  \[ \uExt^m_X(\wC, \wD) \ce \bigoplus_{v \in \ZZ} \Ext^m_X(\wC, \wD(v)) = \bigoplus_{v \in \ZZ} H^m\RHom_X(\wC, \wD(v)). \]
  Finally, we write $\ushHom_X(\wC, \wD)$ for the complex of sheaves associated to $\uHom_R(C, D)$.   
\end{conventions}

\section{Computing global $\Ext$}\label{sec:proof}

The following result of Smith provides the theoretical foundation for the algorithm used by \verb|Macaulay2| for computing $\Ext$ between coherent sheaves on projective varieties, and in particular for computing sheaf cohomology over such varieties.

\begin{theorem}[{\cite[Theorem 1]{Smith00}}]\label{thm:smith}
  Let $M$ and $N$ be finitely generated graded $R$-modules, $m$ an integer, and $\ell \ce \min\{\dim(N), m\}$. If $r \ge \max \{\oa^S_i(N) - i \text{ $\on{:}$ } n - \ell \le i \le \pd_S(N)\} - m + 1$, then there is an isomorphism
  $ \bigoplus_{v \ge 0} \Ext^m_X(\wM, \wN(v)) \cong \uExt^m_R(M_{\ge r}, N)_{\ge 0}$
  of graded $R$-modules.
\end{theorem}

Our goal is to generalize Theorem~\ref{thm:smith} to the context of complexes of coherent sheaves on a projective variety. We will work under the following setup:

\begin{setup}\label{setup}
  Let $C$ and $D$ be bounded complexes of finitely generated $R$-modules and $F$ the minimal $R$-free resolution of $C$, which exists and is unique up to isomorphism by~\cite[Proposition 4.4.1]{Roberts98}. 
  Fix $m \in \ZZ$ and set $\ell \ce \min\{\dim(D) + \sup(D), m + \sup(C) - \inf(D)\}$.
\end{setup}

\begin{remark}
Given integers $d_0, \dots, d_n > 0$, let $\PP(d_0, \dots, d_n)$ denote the stack quotient of~$\mathbb{A}^{n+1}~\setminus~\{0\}$ by the action of the multiplicative group $k^*$ given by $t \cdot (a_0, \dots, a_n) = (t^{d_0}a_0, \dots, t^{d_n}a_n)$. All of the results in this section can be generalized to the setting where the degree of each variable in $S = \k[x_0, \dots, x_n]$ is $\deg(x_i) = d_i$, the ideal $I \subseteq S$ is homogeneous with respect to this grading, and $X$ is the substack of $\PP(d_0, \dots, d_n)$ determined by $I$. We restrict ourselves in this paper to the setting of subvarieties of projective space for the sake of simplicity. 
\end{remark}


We begin by recording some standard technical results on derived global sections of complexes of sheaves. We refer the reader to \cite[Section 3.3]{Huybrechts06} for an introduction to this topic. 

\begin{lemma}\label{lem:dirim}
Let $G$ be a bounded complex of finitely generated $R$-modules. We have $\RR^m\Gamma(X,\wG) = 0$ for $m \ge \dim(G) + \sup(G)$.
\end{lemma}

\begin{proof}
  There is a convergent hypercohomology spectral sequence $E_1^{p, q} = H^q(X, \wG^p) \Rightarrow \RR^{p+q}\Gamma(X, \wG)$. Since $H^q(X, \wG^p) = 0$ for $q \ge \dim(G)$, the conclusion follows.
\end{proof}

\begin{lemma}\label{lem:ext}
  Let $\cG$ and $\cH$ be complexes of $\OO_X$-modules and $\cF$ a locally free resolution of $\cG$.
  There is a canonical quasi-isomorphism $ \RHom_X(\cG, \cH) \simeq \RR\Gamma(X, \ushHom_X(\cF, \cH)). $
\end{lemma}
\begin{proof}
  We have a quasi-isomorphism $\ushHom_X(\cF, \cH) \simeq \RR\ushHom_X(\cG, \cH)$~\cite[pages 76-77]{Huybrechts06}, and there is an adjunction quasi-isomorphism $\RR\ushHom_X(\cG, \cH) \simeq \RR\ushHom_X(\OO_X, \RR\ushHom_X(\cG, \cH))$. Passing to derived global sections finishes the proof.
\end{proof}

We arrive at a canonical map
\begin{equation}\label{eqn:gammamap}
  \Gamma(X, \ushHom_X(\cF, \cH)) \to \RHom_X(\cG, \cH)
\end{equation}
given by composing the natural map $\Gamma(X, \ushHom_X(\cF, \cH)) \to \RR\Gamma(X, \ushHom_X(\cF, \cH))$ with the quasi-isomorphism in Lemma~\ref{lem:ext}.

We now establish several technical statements analogous to results used in Smith's proof of Theorem~\ref{thm:smith}. We begin by restating a result from \cite{Smith00}:

\begin{lemma}[\cite{Smith00} Lemma 2.1]\label{lem:CohomologyVanishingByTwist}
  Let $N$ be a finitely generated graded $R$-module and $m > 0$. If $ v \ge \oa^S_{n-m}(N)-n$, then $ H^m(X,\wN(v)) = 0. $
\end{lemma}

We next extend \cite[Proposition 2.4]{Smith00}. Recall that we use the notation from Setup~\ref{setup}.

\begin{proposition}\label{prop:SheafExtToHom}
  Given $u \in \ZZ$, let $\lambda_u \ce \max\{m - u -\sup(D), -\sup(C)\}$.  Suppose $v \in \ZZ$ satisfies both of the following inequalities:
  \begin{align*}
    v \ge & \max\left\{ \; \oa^S_{n-u}(D^{m - u - s}) \; - \, \ua^R_s(C) \; \text{ $\on{:}$ } 1 \le u \le \ell, \text{ } \lambda_u \le s \le m - u -\inf(D) \right\} - n, \\
    v \ge & \max\left\{\oa^S_{n-u + 1}(D^{m - u - s}) - \ua^R_s(C) \text{ $\on{:}$ } 2 \le u \le \ell, \text{ } \lambda_u \le s \le m - u -\inf(D) \right\} - n.
  \end{align*}
  The map $\Gamma(X, \ushHom_X(\wF, \wD(v))) \to \RHom_X(\wC, \wD(v))$ in \eqref{eqn:gammamap} induces an isomorphism
  \[
  H^{m}\left(\Gamma(X,\ushHom_X(\wF, \wD(v)))\right)
  \cong \Ext^m_X(\wC, \wD(v)).
  \]
\end{proposition}

\begin{proof}
  Set $\cE \ce \ushHom_X(\wF,\wD(v))$. By Lemma~\ref{lem:ext}, it suffices to show that the canonical map $f \co H^m\left(\Gamma(X,\cE )\right) \to \Ext^m_X(\wC,\wD(v))$ is an isomorphism.
  Let $H^q(X,\cE)$ denote the \textit{complex} of $\k$-vector spaces obtained by applying the functor $H^q(X, -)$ to each term of the complex~$\cE$. Writing  $F^{-i}  = \bigoplus_{j \in \ZZ} R(-j)^{\beta^R_{i, j}(C)}$, we have
  \begin{equation}\label{eqn:E2}
    H^q(X, \cE)^t = \bigoplus_{t' + t'' = t} \bigoplus_{j \in \ZZ} H^q(X, \wD^{t''}(j + v)^{\beta^R_{t', j}(C)}).
  \end{equation}
  There is a convergent spectral sequence
  \begin{equation*}
    E_2^{p,q}= H^{p}(H^q(X,\cE))\Rightarrow \RR^{p+q}\Gamma(X,\cE).
  \end{equation*}
  By Lemma~\ref{lem:ext}, it suffices to show that $\bigoplus_{p+q = m} E_{\infty}^{p,q} \cong E_2^{m,0}$; indeed, the resulting isomorphism $E_2^{m,0} = E_{\infty}^{m,0} \xra{\cong} \RR^m\Gamma(X, \cE)$ coincides with the canonical map $f$. Since $E_2^{p,q} = 0$ for $q < 0$ or $p <  \inf(D)-\sup(C)$, we need only show $E_2^{m - u, u} = 0$ for $1 \le u \le m +\sup(C) - \inf(D)$ and $E_2^{m - u, u-1} = 0$ for $2 \le u \le m +\sup(C) - \inf(D)$. Using~\eqref{eqn:E2}, it is thus sufficient to show the following two conditions hold for all $j \in \ZZ$:
  \begin{align*}
    H^u(X,\, \wD^{t''}(j + v)^{\beta^R_{t', j}(C)}) &= 0,
    \qquad 1 \le u \le m +\sup(C) - \inf(D),\quad t' + t'' = m-u; \\
    H^{u-1}(X,\, \wD^{t''}(j + v)^{\beta^R_{t', j}(C)}) &= 0,
    \qquad 2 \le u \le m +\sup(C) - \inf(D),\quad t' + t'' = m-u.
  \end{align*}
  By Lemma~\ref{lem:dirim}, the cohomology groups in each of these two conditions vanish in cohomological degree at least $\dim(D)+\sup(D)$. Recalling from Setup~\ref{setup} that $\ell \ce \min\{\dim(D) + \sup(D), m +\sup(C) - \inf(D)\}$, we need only show the following refinements of the above two conditions hold for all $j \in \ZZ$:
    \begin{align*}
      H^u(X,\, \wD^{t''}(j + v)^{\beta^R_{t', j}(C)}) &= 0,
      \qquad 1 \le u \le \ell ,\quad \lambda_u \le t' \le m-u - \inf(D); \\
      H^{u-1}(X,\, \wD^{t''}(j + v)^{\beta^R_{t', j}(C)}) &= 0,
      \qquad 2 \le u \le \ell,\quad \lambda_u \le t' \le m-u - \inf(D).
    \end{align*}
Lemma~\ref{lem:CohomologyVanishingByTwist} and our bounds on $v$ therefore imply the result. 
\end{proof}

Let $\m$ denote the homogeneous maximal ideal of $S$. For $j \in \ZZ$, let $H^j_\m(\uHom_R(F, D))$ denote the complex of graded $R$-modules obtained by applying the local cohomology functor $H^j_\m(-)$ to the terms of $\uHom_R(F, D)$. We have an exact sequence of complexes of graded $R$-modules:

\begin{equation}\label{eqn:4term}
  0 \to H^0_\m(\uHom_R(F, D)) \to \uHom_R(F, D) \to \bigoplus_{v \in \ZZ} \Gamma(X,\ushHom_X(\wF,\wD(v))) \to H^1_\m(\uHom_R(F, D)) \to 0.
\end{equation}
The following is a generalization of~\cite[Proposition 2.5]{Smith00}:

\begin{proposition}\label{prop:LD}
  Let $\gamma \ce \min\{\sup(D), \sup(C) + m \}$. If
  $$ e \ge \max\{ \max\{\oa^S_n(D^i), \oa^S_{n+1}(D^i)\} - \ua^R_{m-i}(C) \text{ $\on{:}$ } \inf(D) \le i \le \gamma\} - n, $$
  then the canonical morphism
  $ \uHom_R(F, D) \to \bigoplus_{v \in \ZZ} \Gamma(X,\ushHom_X(\wF,\wD(v))) $
  of complexes of graded $R$-modules induces an isomorphism of graded $R$-modules
  $$ \uExt_R^m(C, D)_{\ge e} \cong \bigoplus_{v \ge e} H^{m}\left(\Gamma(X,\ushHom_X(\wF,\wD(v)))\right). $$
\end{proposition}
\begin{proof}
  Fix $t \in \{0, 1\}$. Using the exact sequence~\eqref{eqn:4term}, it suffices to show that the cohomological degree $m$ term of $H^t_\m(\uHom_R(F, D))$ vanishes in internal degrees less than $e$. 
 As in the proof of Proposition~\ref{prop:SheafExtToHom}, write $F^{-i} = \bigoplus_{j = 1}^{b_i} R(-j)^{\beta^R_{i, j}(C)}$, where $b_i \ce \rank(F^{-i})$. We have $\uHom_R(F, D)^m \cong \bigoplus_{i = \inf(D)}^{\gamma} \bigoplus_{j = 1}^{b_{m-i}} D^i(j)^{\beta^R_{m-i,j}(C)}$.
  A local duality argument exactly as in the proof of \cite[Proposition 2.5]{Smith00} implies that, if $H_\m^t(D^i)_c \ne 0$, then we have the inequality $c \le \max\{\oa^S_n(D^i), \oa^S_{n+1}(D^i)\} - n -1$. Thus, if
  $H^t_\m(D^i(j))_c \ne 0$, then
  $$ c \le \max\{\oa^S_n(D^i), \oa^S_{n+1}(D^i)\} - \beta^R_{m-i, j}(C) - n - 1 \le \max\{\oa^S_n(D^i), \oa^S_{n+1}(D^i)\} - \ua^R_{m-i}(C) - n-1, $$
  and the result follows.
\end{proof}

Combining the map from \eqref{eqn:gammamap} with the middle map of \eqref{eqn:4term}, we obtain a canonical morphism of complexes of graded $R$-modules
\begin{equation}\label{eqn:keymap}
  \RuHom_R(C, D) \to \bigoplus_{v \in \ZZ} \RHom_X(\wC, \wD(v)).
\end{equation}
Of course, we may also replace $C$ with $C_{\ge r}$ for any $r \in \ZZ$ to obtain canonical maps
\begin{equation}\label{eqn:keyrmap}
  \RuHom_R(C_{\ge r}, D) \to \bigoplus_{v \in \ZZ} \RHom_X(\wC, \wD(v)).
\end{equation}
The following corollary of Propositions~\ref{prop:SheafExtToHom} and~\ref{prop:LD} is a generalization of~\cite[Corollary 2.6]{Smith00}:
\begin{corollary}\label{cor:smithcpx}
  Let $\lambda_u$ be as in Proposition~\ref{prop:SheafExtToHom}, and suppose $e \in \ZZ$ satisfies both of the following:
  \begin{align*}
    e \ge & \max\{ \; \oa^S_{n-u}(D^{m-u-s}) \; - \, \ua^R_s(C) \; \text{ $\on{:}$ } 0 \le u \le \ell, \quad\quad\quad\, \text{ } \lambda_u \le s \le m - u -\inf(D) \} - n, \\
    e \ge & \max\{\oa^S_{n-u + 1}(D^{m-u-s}) - \ua^R_s(C) \text{ $\on{:}$ } 0 \le u \le \ell, \text{ } u \ne 1, \text{ } \lambda_u \le s \le m - u -\inf(D) \}-n.
  \end{align*}
  The map~\eqref{eqn:keymap} induces an isomorphism $\uExt_R^m(C, D)_{\ge e} \cong \bigoplus_{v \ge e} \Ext_X^m(\wC, \wD(v))$.
\end{corollary}
\begin{proof}
  Applying the change of variables $i = m-s$ to the inequalities in Proposition~\ref{prop:LD} gives the $u = 0$ cases of the inequalities in the statement. Propositions~\ref{prop:SheafExtToHom} and~\ref{prop:LD} thus imply the result.
\end{proof}

We are now ready to prove the main result of this section:

\begin{theorem}\label{thm:smithcpx}
  Suppose the following inequality holds:
  $$ r \ge \max\{\oa^S_{i}(D^j) \text{ $\on{:}$ } n-\ell \le i \le \pd_S(D^j), \text{ } \inf(D) \le j \le \min\{ \sup(D) , \sup(C) - n +m+i\} \} - n. $$
  The map~\eqref{eqn:keyrmap} induces an isomorphism
  $ \uExt^m_R(C_{\ge r}, D)_{\ge 0}  \cong\bigoplus_{v \ge 0} \Ext^m_X(\wC, \wD(v)) $
  of graded $R$-modules. Moreover, when $H^j(C) = 0$ for $j \ne 0$, this isomorphism is implied by the combination of the following inequalities:
  \begin{footnotesize}
    \begin{align*}
      r \ge & \max\{\oa^S_{n-\ell}(D^j)+j \; \text{ $\on{:}$ } \; n-\ell \le \pd_S(D^j), \quad\quad\quad\,\, \text{ } \inf(D)  \le j \le \min\{\sup(D), m - \ell \} \} -n + \ell - m, \\
      r \ge & \max\{\oa^S_{i}(D^{j}) +j - i \text{ $\on{:}$ } n-\ell + 1 \le i \le \pd_S(D^j),  \text{ } \inf(D) \le j \le \min\{\sup(D),m - n + i - 1 \}\} + 1 - m.
    \end{align*}
  \end{footnotesize}
\end{theorem}

\begin{remark}
  The inequality in the first statement of Theorem~\ref{thm:smithcpx} implies both inequalities in the second statement.
\end{remark}
\begin{proof}
  Let $\lambda_u$ be as in Proposition~\ref{prop:SheafExtToHom}. By Corollary~\ref{cor:smithcpx}, it suffices to show  the following two inequalities hold:
  \begin{align*}
    0 \ge & \max\{ \; \oa^S_{n-u}(D^{m-u-s}) \; - \, \ua^R_s(C_{\ge r}) \, \text{ : } 0 \le u \le \ell, \quad\quad\quad\, \text{ } \lambda_u \le s \le m - u -\inf(D) \}-n, \\
    0 \ge & \max\{\oa^S_{n-u + 1}(D^{m-u-s}) - \ua^R_s(C_{\ge r}) \text{ : } 0 \le u \le \ell, \text{ } u \ne 1, \text{ } \lambda_u \le s \le m - u - \inf(D) \}-n.
  \end{align*}
  Note that $\ua^R_t(C_{\ge r}) \ge r$ for all $t$.\footnote{When $H^j(C) = 0$ for $j \ne 0$, this inequality can be strengthened to $\ua^R_t(C_{\ge r}) \ge r + t$ for all $t$; we use this observation to prove the second statement of the theorem. Our use of the weaker inequality $\ua^R_t(C_{\ge r}) \ge r$ is the reason why the first statement of our theorem doesn't generalize Smith's Theorem (Theorem~\ref{thm:smith}).}  Applying the change of variables $i \ce n-u$ and $j \ce m - u - s$, one sees that the first inequality is implied by
  $$ r \ge  \max\{\oa^S_{i}(D^j) \text{ : } n-\ell \le i \le n, \text{ } \inf(D)  \le j \le \min\{ \sup(D) , \sup(C) - n +m+i\} \} - n. $$
  Using the change of variables $i \ce n-u+1$ and $j \ce m-u-s$, the second inequality is implied by
  \begin{footnotesize}
    $$ r \ge  \max\{\oa^S_{i}(D^{j}) \text{ : } n+1-\ell \le i \le n+1, \text{ } i \ne n,\text{ } \inf(D) \le j \le  \min\{\sup(D), \sup(C) -n + m +i - 1 \} \} - n. $$
  \end{footnotesize}
  Since $\oa^S_i(D^j) = -\infty$ when $i > \pd_S(D^j)$, the first statement follows. To prove the second statement, we strengthen the  inequality $\ua^R_t(C_{\ge r}) \ge r$ to $\ua^R_t(C_{\ge r}) \ge r + t$  and proceed exactly as above.
\end{proof}

Theorem~\ref{thm:smithcpx} takes on a more pleasing form if $\sup(C)=\sup(D) = 0$; of course one can always shift $C$ and $D$ to ensure this. We state this special case as a corollary, which implies Theorem~\ref{thm:smithcpxintro}.

\begin{corollary}\label{cor:positive}
  Suppose $\sup(C) = \sup(D) = 0$, and assume the following inequality holds:
  $$ r \ge \max\{\oa^S_{i}(D^{j}) \text{ $\on{:}$ } n-\ell \le i \le \pd_S(D^j), \text{ } \inf(D) \le j \le  0 \} - n. $$
  The map~\eqref{eqn:keyrmap} induces an isomorphism
  $ \uExt^m_R(C_{\ge r}, D)_{\ge 0} \cong \bigoplus_{v \ge 0} \Ext^m_X(\wC, \wD(v)) $
  of graded $R$-modules. Moreover, when $H^j(C) = 0$ for $j \ne 0$, we need only check the following inequalities:
  \begin{align*}
    r \ge & \max\{\oa^S_{n-\ell}(D^j)+j \text{ $\on{:}$ } n-\ell \le \pd_S(D^j), \text{ } \inf(D) \le j \le  0 \} -n + \ell - m, \\
    r \ge & \max\{\oa^S_{i}(D^{j}) +j- i  \text{ $\on{:}$ } n-\ell + 1 \le i \le \pd_S(D^j),  \text{ } \inf(D) \le j \le  0 \}+1 - m.
  \end{align*}
\end{corollary}

\begin{remark}
  When $C$ and $D$ are concentrated in cohomological degree zero, Corollary~\ref{cor:positive} also recovers (and is in fact slightly stronger than) Theorem~\ref{thm:smith}; this special case of Theorem~\ref{thm:smithcpx} also easily follows from Smith's proof of Theorem~\ref{thm:smith}.
\end{remark}

The following algorithm is adapted from Theorem~\ref{thm:smithcpxintro} and implemented in \verb|Macaulay2|.

\begin{algorithm}[H]
  \caption{\!\!\textbf{.} \texttt{RHom(F,G)} 
  }\label{alg:RHom}
  \begin{algorithmic}[1]
    \smallskip
    \Require two bounded complexes of finitely generated graded modules $C$ and $D$ over $R = S/I$ such that $\sup(C) = \sup(D) = 0$, representing complexes $\cF$ and $\cG$ of coherent sheaves on $X = \Proj R$.
             
    \Ensure  a complex of $\k$-vector spaces whose $m^{\th}$ cohomology is $\Ext^m_X(\cF, \cG)$.
    \State $r \ce -\infty$
    \State $n \ce \dim S - 1$ 
    \ForAll { $\inf(D) \le j \le 0$ }
      \State $F \ce$ minimal $S$-free resolution of $D^j$
      \ForAll { $0 \le i \le \text{length}(F)$ }
        \State $d \ce $ highest degree of a generator of $F_i$
        \State $r \,= \max\{r \,, \; d - n \}$
      \EndFor
    \EndFor
    \State $C' \ce$ minimal free resolution of the truncation $C_{\ge r}$
    \State \Return $\uHom_R(C', D)_0$
  \end{algorithmic}
\end{algorithm}


\section{Examples}\label{sec:examples}

We give two examples of $\Ext$ calculations for complexes using Theorem~\ref{thm:smithcpxintro}.





\begin{example}\label{ex:koszul}
Suppose $R = S = \k[x_0, x_1, x_2]$, let $K$ be the Koszul complex on the (non-regular) sequence $x_0^2, x_0x_1$, and let $\cK$ be the complex of sheaves on $\PP^2$ associated to $K$. The complex~$\cK$ has nonzero cohomology in degrees $0$ and $-1$. Let us compute the derived global sections of~$\cK$. As in the proof of Lemma~\ref{lem:dirim}, there is a convergent hypercohomology spectral sequence $E_1^{p,q} = H^q(\PP^2, \cK^p) \Rightarrow \RR^{p+q}\Gamma(\PP^2, \cK)$. We have $E_1^{p,q} = \k$ when $p = q = 0$, $E_1^{p,q} = \k^3$ when $p = -2$ and $q = 2$, and $E_1^{p,q} = 0$ otherwise.
  Thus, the spectral sequence degenerates at page 1. We conclude that
  $\RR^0\Gamma(\PP^2, \cK) \cong \k^4$, and
  $\RR^m\Gamma(\PP^2, \cK) = 0$ for $m \ne 0$.

Let us now recover this calculation using Theorem~\ref{thm:smithcpxintro}. We have:
  $$
  \max\{\oa^S_{i}(K^{j}) \text{ $\on{:}$ } 0 \le i \le \pd_S(K^j), \text{ } \inf(K) \le j \le  0\} - n = 4 - 2 = 2.
  $$
  Thus, if $r \ge 2$, then Theorem~\ref{thm:smithcpxintro} implies $\RR^m\Gamma(\PP^2, \cK) \cong \Ext^m_S(S_{\ge r}, K)$ for all $m$. We verify this for $m = 0$ in \verb|Macaulay2|:

\begin{footnotesize}
\begin{lstlisting}[language=Macaulay2, frame=single]
`\underline{\tt i1}` : kk = ZZ/32003;
`\underline{\tt i2}` : S = kk[x_0,x_1,x_2];
`\underline{\tt i3}` : C = koszulComplex {x_0^2, x_0*x_1};
`\underline{\tt i4}` : F = res module truncate(2, S);
`\underline{\tt i5}` : part(0, HH^0 Hom(F, C))
`\underline{\tt o5}` = kk^4
\end{lstlisting}
\end{footnotesize}
  Moreover, the bound $r \ge 2$ is sharp in this example: we have $\Ext_S^0(S_{\ge 1}, K) \cong \k \ne \k^4$.
\end{example}

\begin{example}
  Let $S = \k[x_0, x_1, x_2]$, $R$ the hypersurface $S / (x_0x_1)$, $K$ the Koszul complex on $x_0$ over $R$, and $\cK$ its sheafification on $X = \Proj(R) \subseteq \PP^2$. The complex $\cK$ has nonzero cohomology in degrees $0$ and $-1$.
  We now compute $\Ext^*_X(\cK, \cK(1))$. Let $\cC$ denote the complex $\ushHom_X(\cK, \cK(1))$, which has the form
\begin{equation*}
\OO_X \xrightarrow{\begin{pmatrix}x_0\\ x_0\end{pmatrix}} \OO_X(1)^2 \xrightarrow{\begin{pmatrix}x_0&-x_0\end{pmatrix}} \OO_X(2),
\end{equation*}
where $\OO_X$ is in cohomological degree $-1$. By Proposition~\ref{prop:SheafExtToHom}, we have
$
\Ext^*_X(\cK,\cK(1))\cong \RR^*\Gamma(X,\mathcal{C}).
$
Each term of $\cC$ has sheaf cohomology only in degree 0, and so a spectral sequence argument implies that its derived global sections are given by the cohomology of the induced complex
  \begin{equation}
  \label{eqn:short}
  0\to H^0(X, \OO_X) \to H^0(X,\OO_X(1))^2 \to H^0(X,\OO_X(2)) \to 0.
  \end{equation}
  A direct calculation thus yields that $\Ext_X^m(\cK, \cK(1)) \cong \k^3$ when $m \in \{0,1\}$ and $\Ext_X^m(\cK, \cK(1)) = 0$ otherwise. Let us now confirm this via Theorem~\ref{thm:smithcpxintro}. In this case, Theorem~\ref{thm:smithcpxintro} states that we may take $r = 0$ to be our truncation degree; that is, we need not truncate at all. Theorem~\ref{thm:smithcpxintro} therefore says that $\Ext^m_X(\cK, \cK(1)) \cong \Ext^m_R(K, K(1))$, which recovers our observation that the cohomology of the complex~\eqref{eqn:short} computes the desired Ext groups. Finally, we carry out this calculation in \verb|Macaulay2|. In the final output, note that homological indexing is the default for \verb|Macaulay2|.
  \begin{footnotesize}
\begin{lstlisting}[language=Macaulay2, frame=single]
`\underline{\tt i6}` : R = S/(x_0*x_1);
`\underline{\tt i7}` : K = koszulComplex {x_0};
`\underline{\tt i8}` : part(0, HH Hom(K, K ** R^{1}))
`\underline{\tt o8}` = kk^3 <- kk^3 <- 0
\end{lstlisting}
\end{footnotesize}
\end{example}




\section{Exceptional collections of complexes}\label{sec:exc-cols}

\let\<\langle
\let\>\rangle

An object $\cC$ in the bounded derived category $\Db(X)$ is called \emph{exceptional} if $\Ext^m_X(\cC, \cC)$ is zero for $m \neq 0$ and is isomorphic to $\k$ for $m = 0$. A sequence of exceptional objects $\cC_1, \dots, \cC_t$ in $\Db(X)$ is called an \emph{exceptional collection} if $\Ext^m_X(\cC_i, \cC_j) = 0$ when $i > j$ for all $m$, and a \emph{strong} exceptional collection if further $\Ext^m_X(\cC_i, \cC_j) = 0$ when $i < j$ and $m \neq 0$. Exceptional collections are an important topic of study in algebraic geometry; see e.g. \cite{Huybrechts06} for background. In the following examples we highlight applications of Theorem~\ref{thm:smithcpxintro} in working with exceptional collections consisting of complexes of coherent sheaves in $\Db(X)$.

\newcommand{\ev}{\mathrm{ev}}
\begin{example}[Mutations of exceptional collections]
  Given an exceptional pair $(\cE,\cF)$, define the left mutation $L_\cE(\cF)$ and right mutation $R_\cF(\cE)$ to be objects fitting in the distinguished triangles
  \begin{align*}
    L_\cE(\cF) \to \RHom_X(\cE, \cF) \otimes \cE \xra{\ev} \cF \quad\text{and}\quad
    \cE \xra{\ev^*} \RHom_X(\cE, \cF)^* \otimes \cF \to R_\cF(\cE),
  \end{align*}
  where $\ev$ (resp.~$\ev^*$) is the derived (co)evaluation map and $\RHom_X(\cE,\cF)^*$ is the dual complex of $\kk$-vector spaces. In particular, $(L_\cE(\cF), \cE)$ and $(\cF, R_\cF(\cE))$ are again exceptional pairs \cite[\S2]{Bondal90}. Mutations induce an action from the braid group of $t$ strings on exceptional collections with $t$ objects. The derived evaluation map $\ev$ is given by a direct sum of maps $\cE[m]\to\cF$ for each basis vector in $H^m\RHom_X(\cE,\cF)$, while in $\ev^*$ the corresponding direct summand is $\cE\to\cF[-m]$. Thus both maps are explicitly computable with the help of Algorithm~\ref{alg:RHom}.

  As a basic application, we compute the exceptional collection which is left dual to Beilinson's exceptional collection $\OO_{X}, \OO_{X}(1), \OO_{X}(2)$ on $X = \PP^2$ via three left mutations as follows:
  \newcolumntype{L}{>{$}l<{$}} 
  \begin{center}
  \renewcommand{\arraystretch}{1.5}
  \begin{tabular}{LLLL} \hline
    \OO_{X} & \OO_{X}(1) & \OO_{X}(2) \\ \hline
    \noalign{\vskip 1mm}
    \OO_{X} & L_{\OO_{X}(1)}\OO_{X}(2) = {\OO_{X}(1)}^3\to\OO_{X}(2) & \OO_{X}(1) \\ \hline
    \noalign{\vskip 1mm}
    L_{\OO_{X}}L_{\OO_{X}(1)}\OO_{X}(2) = {\OO_{X}}^3\to{\OO_{X}(1)}^3\to\OO_{X}(2) & \OO_{X} & \OO_{X}(1) \\ \hline
    \noalign{\vskip 1mm}
    L_{\OO_{X}}L_{\OO_{X}(1)}\OO_{X}(2) = {\OO_{X}}^3\to{\OO_{X}(1)}^3\to\OO_{X}(2) & L_{\OO_{X}}\OO_{X}(1) = {\OO_{X}}^3\to\OO_{X}(1) & \OO_{X} \\ \hline
  \end{tabular}
  \end{center}
  On $\PP^n$\!, this recovers the full strong exceptional collection $\Omega^n(n), \dots, \Omega^1(1), \OO \in \Db(\PP^n)$ \cite{Beilinson78}. Observe that the exterior powers of the cotangent sheaf are given by truncations of the Koszul complex in two different ways which may be achieved via left and right mutations:
  \begin{align*}
    0 \gets \Omega^a \gets \bigwedge^{a+1} \left[\OO^{n+1}(-1)\right] \gets \bigwedge^{a+2} \left[\OO^{n+1}(-1)\right] &\gets \dots \gets \OO(-n-1) \gets 0 \\
    0 \gets \OO^1 \gets \dots \gets \bigwedge^{a-1} \left[\OO^{n+1}(-1)\right] \gets \bigwedge^a \left[\OO^{n+1}(-1)\right] &\gets \Omega^a \gets 0.
  \end{align*}

  In the example above, we can display the ranks of the extension groups among the collection of complexes on the last row as
  \[
  \begin{pmatrix}
    1&3&3\\
    0&1&3\\
    0&0&1
  \end{pmatrix},
  \]
  where each term $c_m T^m$ in the matrix entry $(i,j)$ indicates that the rank of $\Ext_X^m(\cE_i,\cE_j)$ is $c_m$ (in general, entries of such an extension table can be Laurent polynomials in $T$).
  Hence the Beilinson collection~$\{ \Omega^2(2), \Omega(1), \OO \}$ forms a strong exceptional collection on $\PP^2$.
\end{example}

\begin{example}[Spherical twists and exceptional collections]
  Assume $X$ is a smooth projective variety. An object $\cE\in\Db(X)$ is called \emph{spherical} if $\cE \otimes \omega_X \cong \cE$ and $\Hom_X(\cE, \cE[i])$ is $\k$ when $i=0$ or $i=\dim(X)$ and zero otherwise. Such an object induces an auto-equivalence of the derived category
  \begin{align*}
    T_\cE\colon \Db(X) &\xra{\simeq} \Db(X), \\
    \text{given by } \cF &\mapsto L_\cE(\cF)[1]
  \end{align*}
  known as a \emph{spherical twist} \cite[Prop.~8.6]{Huybrechts06}. (We are abusing notation slightly, as $L_\cE(\cF)$ was only defined above in the case where $(\cE, \cF)$ is an exceptional collection.) In particular, applying the spherical twist to every term in an exceptional collection yields another exceptional collection.

  As an application, we construct an exceptional collection on $\FF_2$, the Hirzebruch surface of type 2,\linebreak consisting of objects not quasi-isomorphic to sheaves. Suppose $E$ and $H$ are the exceptional divisor on $\FF_2$ and the pullback of the hyperplane divisor, respectively, and consider the coordinate ring of the embedding $X = \Proj R \xhookrightarrow{} \PP^5$ corresponding to the ample line bundle $\OO_{\FF_2}(E+3H)$
  \[ R = \k[x_0,\dots,x_5]/(x_{4}^{2}-x_{3}x_{5}, x_{3}x_{4}-x_{2}x_{5}, x_{1}x_{4}-x_{0}x_{5}, x_{3}^{2}-x_{2}x_{4}, x_{1}x_{3}-x_{0}x_{4}, x_{1}x_{2}-x_{0}x_{3}). \]
  Following \cite{IshiiOkawaUehara21}, the line bundle $\OO_E$ on the $(-2)$-curve in $\FF_2$ is a spherical object.
  The pushforwards of terms of an exceptional collection on $\FF_2$ to $X$ via this embedding are
  \newcolumntype{C}{>{$}c<{$}} 
  \begin{center}
  \renewcommand{\arraystretch}{1.5}
  \begin{tabular}{CCCC}\hline
    \OO_{\FF_2} & \OO_{\FF_2}(H) & \OO_{\FF_2}(2H+E) & \OO_{\FF_2}(3H+E) \\
    \hline
    \noalign{\vskip 1mm}
    A = \OO_X &
    B = \mathrm{image}\Big({\OO_X}^2 \xrightarrow{\;m^T} {\OO_X(1)}^4\Big) &
    C = \mathrm{image}\Big({\OO_X}^4 \xrightarrow{\;m\;} {\OO_X(1)}^2\Big) &
    D = \OO_X(1) \\[1ex]
    \hline
  \end{tabular}
  \end{center}
  where the matrix $m$ in differentials of $B$ and $C$ is
  \[ m = \begin{pmatrix} x_{5}&x_{4}&x_{3}&x_{1}\\ x_{4}&x_{3}&x_{2}&x_{0} \end{pmatrix}. \]
  In this embedding, one can verify that $\OO_E = \mathrm{coker}(x_2 \; x_3 \; x_4 \; x_5)$ is the spherical object. In order to compute the spherical twists, we need to compute $\RHom_X(\OO_E, -)$ for each term of the collection. The only nontrivial extensions are $\Ext_X^2(\OO_E, A)$ and $\Ext_X^2(\OO_E, C)$. Therefore the spherical twists, whose cohomology is not always concentrated in a single degree, are given by
  \begin{center}
  \renewcommand{\arraystretch}{1.5}
  \begin{tabular}{CCCC}\hline
    T_{\OO_E}(A) & T_{\OO_E}(B) & T_{\OO_E}(C) & T_{\OO_E}(D) \\
    \hline
    \noalign{\vskip 1mm}
    \cone\left(\Ext^2(\OO_E, A) \otimes \OO_E[-2] \to A \right) &
    B &
    \cone\left(\Ext^2(\OO_E, C) \otimes \OO_E[-2] \to C \right) &
    D \\[1ex]
    \hline
  \end{tabular}
  \end{center}
  We have $H^*(T_{\OO_E}(A)) = A \oplus \OO_E[-1]$, and $H^*(T_{\OO_E}(C)) = C \oplus \OO_E[-1]$. We can verify that this is again a strong exceptional collection by computing its extension table:
  \[
  \begin{pmatrix}
    1&2&4&6\\
    0&1&2&4\\
    0&0&1&2\\
    0&0&0&1
  \end{pmatrix}.
  \]
\end{example}




\printbibliography

\Addresses

\end{document}